\documentclass[12pt,A4paper]{article}
\usepackage{amsfonts}
\usepackage{mathrsfs}
\usepackage{amssymb,amsfonts,amsmath}
\usepackage{epsfig}
 \usepackage[rotateright]{rotating}
\usepackage{ifpdf}
\newtheorem{thm}{Theorem}[section]
\newtheorem{cor}[thm]{Corollary}

\newtheorem{defn}[thm]{Definition}

\makeatletter \@addtoreset{equation}{section}

\def\E{{\rm E}}\def\Var{{\rm Var}}

\def\pf{\noindent {\it Proof.\ }}
\def\qed{\nopagebreak\hfill{\rule{4pt}{7pt}}}
\def\Des{\mathrm{Des}}
\def\Asc{\mathrm{Asc}}
\def\asc{\mathrm{asc}}
\def\des{\mathrm{des}}

\def\Exc{\mathrm{Exc}}
\def\exc{\mathrm{exc}}
\def\wexc{\mathrm{wexc}}

\def\sg{\sigma}

\begin{document}
\begin{center}
{\large\bf  Derangement Polynomials and Excedances of Type $B$}\\
\end{center}

\begin{center}
William Y. C. Chen$^{1}$, Robert L. Tang$^2$ and
Alina F. Y. Zhao$^{3}$\\[6pt]
Center for Combinatorics, LPMC-TJKLC\\
Nankai University, Tianjin 300071, P. R. China\\[5pt]
$^{1}${chen@nankai.edu.cn}, $^{2}${tangling@cfc.nankai.edu.cn},
 $^{3}${zfeiyan@cfc.nankai.edu.cn}
\end{center}

\vskip 3mm \noindent \textbf{Abstract.} Adopting the definition of
excedances of type $B$ due to Brenti, we give a type $B$ analogue of
the $q$-derangement polynomials. The connection between
$q$-derangement polynomials and Eulerian polynomials naturally
extends to the type $B$ case. Based on this relation,  we derive
some basic properties of the $q$-derangement polynomials of type
$B$, including the generating function formula, the Sturm sequence
property, and the asymptotic normal distribution. We also show that
the $q$-derangement polynomials are almost symmetric in the sense
that the coefficients possess the spiral property.



\section{Introduction}

In this paper, we define a type $B$ analogue of the $q$-derangement
polynomials introduced by Brenti  \cite{Bre90} by $q$-counting
derangements with respect to
 the number of excedances of type $B$, also introduced by Brenti
 \cite{Bre94}. We give some basic properties of these polynomials.
 It turns out that the connection
 between the $q$-derangement polynomials and the Eulerian polynomials
 naturally extends to the type $B$ case, where the type
 $B$ analogue of Eulerian polynomial has been given by Brenti \cite{Bre94},
 and has been further studied by Chow and Gessel in \cite{Chow2}.

Let us now recall some definitions. Let $\mathfrak{S}_n$ be the set
of permutations of $[n]=\{1, 2, \ldots , n\}$. For each $\sg\in
\mathfrak{S}_n$, the \emph{descent} set and the \emph{excedance} set
of $\sg=\sg_1\sg_2\cdots \sg_n$ are defined respectively as follows,
\begin{align*}
\Des(\sg)&= \{i\in [n-1] : \sg_i>\sg_{i+1}\},\\[6pt]
\Exc (\sg)&=\{i \in [n-1] : \sg_i > i \}.
\end{align*}
The descent number and excedance number are defined respectively by
\[ \des(\sg)= |\Des(\sg)|, \qquad \exc(\sg)= |\Exc(\sg)|.\] The
\emph{Eulerian polynomials} \cite{Stanley, Wachs} are defined by
 \[A_n(q)= \sum_{\sg \in \mathfrak{S}_n}
q^{\des(\sg)+1} =\sum_{\sg \in \mathfrak{S}_n} q^{\exc(\sg)+1},
\quad n\geq 0.
\]
The Eulerian polynomials have the following generating function
\begin{equation}\label{gfa}
\sum_{n\geq
0}{A_n(q)}\frac{t^n}{n!}=\frac{(1-q)e^{qt}}{e^{qt}-qe^t}.
\end{equation}

A permutation $\sg=\sg_1\sg_2\cdots \sg_n$ is a \emph{derangement}
if  $\sg_{i}\neq i$ for any $i\in [n]$.
 The set of derangements on $[n]$ is denoted by $D_n$.
 Brenti \cite{Bre90} defined  the $q$-derangement polynomials of type $A$ by
 \[d_n(q)=\sum_{\sg \in D_n}q^{\exc(\sg)},\]
and proved that $d_n(q)$ is symmetric and unimodal for $n\geq 1$.
The following formulas (\ref{e-d-a}) and (\ref{gf}) are derived by
Brenti \cite{Bre90}.

\begin{thm}\label{zhang2} For $n\geq 0$,
\begin{equation}
\label{e-d-a} d_n(q)=\sum_{k=0}^n(-1)^{n-k}{n\choose k}A_k(q).
\end{equation}
\end{thm}

\begin{thm} We have
\begin{equation}\label{gf} \sum_{n\geq
0}d_n(q)\frac{t^n}{n!}=\frac{1}{1-\sum_{n\geq
2}(q+q^2+\cdots+q^{n-1})t^n/n!}.
\end{equation}
\end{thm}

A combinatorial proof of the above identity is given by Kim and Zeng
\cite{Kim} based on a decomposition of derangements. Brenti further
proposed the conjecture that $d_n(q)$ have only real roots for
$n\geq 1$, which has been proved independently  by  Zhang
\cite{Zhang1}, and  Canfield as mentioned in \cite{Bre93}.

\begin{thm}\label{zhang1}
The polynomials $d_n(q)$ form a Sturm sequence. Precisely, $d_n(q)$
has $n$ distinct non-positive real roots, separated by the roots of
$d_{n-1}(q)$.
\end{thm}

The following recurrence relation is  given by Zhang \cite{Zhang1},
which has been used to prove Theorem \ref{zhang1}.

\begin{thm}\label{zhang3}
For $n\geq 2$, we have
$$d_n(q)=(n-1)qd_{n-1}(q)+q(1-q){d_{n-1}^\prime}(q)+(n-1)qd_{n-2}(q).$$
\end{thm}

This paper is motivated by finding the right definition of a type
$B$ analogue of the $q$-derangement polynomials of type $A$ so that
we can get analogous properties to the above theorems for the type
$A$ $q$-derangement polynomials. We discover that the notion of
excedances of type $B$ introduced by Brenti serves as the right
choice for type $B$ derangement polynomials, although there are
several possibilities to define type $B$ excedances,
 see \cite{Bre94,Chow1,ste}.
Nevertheless, it should be
 noted that the type $B$ derangement polynomials are not symmetric
 compared with the case of type $A$. On the other hand, we will be
 able to show that they are almost symmetric in the sense that their
 coefficients have the spiral property.

This paper is organized as follows. In Section 2, we recall Brenti's
definition of type $B$ excedances, and  present the definition of
$q$-derangement polynomials of type $B$, denoted by $d_n^B(q)$. In
Section 3, we establish the connection between
 the derangement polynomials of type $B$ and the Eulerian polynomials
 of type $B$.
This leads to a generating function formula for type $B$ derangement
polynomials. We then extend the $U$-algorithm and $V$-algorithm
given by Kim and Zeng \cite{Kim} to derangements of type $B$. This
gives a combinatorial interpretation of the generating function
formula. In Section 4, we prove that the polynomials $d_n^B(q)$ form
a Sturm sequence. We also show that the coefficients of $d_n^B(q)$
possess the spiral property. In Section 5, by using Lyapunov's
theorem we deduce that the limiting distribution of the coefficients
of $d_n^B(q)$ is normal.

\section{The  Excedances of Type $B$}

In this section, we recall Brenti's definition of type $B$
excedances, and  give the definition of the $q$-derangement
polynomials of type $B$. We adopt the notation and terminology on
permutations of type $B$, or signed permutations, as given in
\cite{Chow1}. Let $B_n$ be the hyperoctahedral group on $[n]$. We
can view the elements of $B_n$ as \emph{signed permutations} of
$[n]$, written as $\sigma=\sigma_1\sigma_2\cdots\sigma_n$, in which
some elements are associated with a minus sign. We also express a
negative element $-i$ in the form  $\bar{i}$.

The \emph{type $B$ descent set} and the \emph{type $B$ ascent set}
of a signed permutation $\sg$ are defined by
\begin{align*}
\Des_B(\sg)&=\{ i \in [0,n-1]:\sg_i>\sg_{i+1}\},\\[6pt]
\Asc_B(\sg)&=\{ i \in [0,n-1]:\sg_i<\sg_{i+1}\},
\end{align*}
where we set $\sg_0=0$. The type $B$ descent and ascent numbers are
given by \[ \des_B(\sg)=|\Des_B(\sg)|,  \qquad
\asc_B(\sg)=|\Asc_B(\sg)|.\]

 A \emph{derangement of type $B$}
on $[n]$ is a signed permutation $\sg=\sg_{1}\sg_{2}\cdots\sg_{n}$
such that $\sg_{i}\neq i$, for all $i\in [n]$. The fixed point is a
position $i$ such that $\sg_i=i$. The set of  derangements in $B_n$
is denoted by $D^B_n$.

 Let us recall the definitions of excedances and
weak excedances of type $B$ introduced by Brenti \cite{Bre94}. For
further information on statistics on signed permutations,
 see \cite{Bre94, Chow2, rei1}.

\begin{defn}
Given $\sg \in B_n$ and $i \in [n]$, we say that $i$ is a type $B$
excedance  of $\sg$ if either $\sg_i=-i$
 or $\sg_{|\sg_i|}>\sg_i$. We denote by $\exc_B(\sg)$
the number of type $B$ excedances  of $\sg$. Similarly,  we say that
$i$ is a type $B$
 weak excedance of $\sg$ if either
$\sg_i=i$ or $\sg_{|\sg_i|}>\sg_i$, and we  denote by
 $\wexc_B(\sg)$ the number of type $B$
weak excedances of $\sg$.
\end{defn}

 Based on the above definition of type
$B$ excedances,
 we  define the type $B$ analogue of the
$q$-derangement polynomials.

\begin{defn} The type $B$ derangement polynomials $d_n^B(q)$
are defined by
\begin{equation}\label{def-d_n}
d_n^B(q)=\sum_{\sg \in
D^B_n}q^{\exc_B(\sg)}=\sum_{k=1}^nd_{n,\,k}\,q^k, \quad n\geq 1,
\end{equation}
where $d_{n,\,k}$ is the number of derangements in $D_n^B$ with
exactly $k$  excedances of type $B$. For $n=0$, we define
$d_0^B(q)=1$.
\end{defn}

Below are the polynomials $d_n^B(q)$ for $n\leq 10$:
\allowdisplaybreaks
\begin{eqnarray*}
d_1^B(q) & = & q,\\[3pt]
d_2^B(q) &  = & 4q+q^2, \\[3pt]
d_3^B(q)& = & 8q+20q^2+q^3, \\[3pt]
 d_4^B(q) & = & 16q+144q^2+72q^3+q^4, \\[3pt]
 d_5^B(q) & = & 32q+752q^2+1312q^3+232q^4+q^5, \\[3pt]
d_6^B(q)& = &64q+3456q^2+14576q^3+9136q^4+716q^5+ q^6, \\[3pt]
d_7^B(q) & = &128q+14912q^2+127584q^3+190864q^4+55624q^5+2172q^6+q^7,\\[3pt]
d_8^B(q) & = &
256q+62208q^2+977920q^3+2879232q^4+2020192q^5\\[3pt]
&&+\,314208q^6+6544q^7+q^8,\\[3pt]
d_9^B(q) & = &
512q+254720q^2+6914816q^3+35832320q^4+49168832q^5\\[3pt]
&&+\,18801824q^6+1697408q^7+19664q^8+q^9,\\[3pt]
d_{10}^B(q) &=&
1024q+1032192q^2+46429440q^3+394153728q^4+937670016q^5\\[3pt]
&&+\,704504832q^6+161032224q^7+8919456q^8+59028q^9+q^{10}.\\[3pt]
\end{eqnarray*}

\section{The Generating Function}

The first result of this section is a formula expressing $d_n^B(q)$
in terms of $B_n(q)$, where $B_n(q)$ are the Eulerian polynomials of
type $B$. This formula is analogous to that of Brenti \cite{Bre90}
for the type $A$ case , and it enables us to derive a formula for
the generating function of $d_n^B(q)$.

 The \emph{Eulerian polynomials of type $B$} are
defined by Brenti \cite{Bre94} based on the number of descents of
type $B$:
\begin{equation}
B_n(q)=\sum_{\sg\in B_n} q^{\des_B(\sg)}.
\end{equation}
Brenti \cite{Bre94} obtained the following formula for the
generating function of the Eulerian polynomials of type $B$, see,
also, Chow and Gessel \cite{Chow2},
\begin{equation}\label{eg}
\sum_{n=0}^\infty
B_n(q)\frac{t^n}{n!}=\frac{(1-q)e^{t(1-q)}}{1-qe^{2t(1-q)}}.
\end{equation}

The following theorem is obtained by  Brenti \cite{Bre94}, which
will be used to establish the formula for $d_n^B(q)$.

\begin{thm}\label{e-we}
There is a bijection $\varphi$: $B_n\rightarrow B_n$ such that
$$\asc_B(\varphi(\sigma))=\wexc_B(\sigma),$$
for any $\sigma\in B_n$.
\end{thm}

The following formula indicates that the notion of  excedances of
type $B$ introduced by Brenti is a right choice for type $B$
derangement polynomials.

\begin{thm}We have
\begin{equation}\label{d-B}
d_n^B(q) = \sum_{k=0}^n (-1)^{n-k}{n \choose k}B_k(q).
\end{equation}
\end{thm}

\pf From Theorem \ref{e-we}, we see that the number of excedances of
type $B$ and the number of descents of type $B$ are equidistributed
on $B_n$. So we deduce that
\begin{equation}
B_n(q)=\sum_{\sg\in B_n}q^{\des_B(\sg)}=\sum_{\sg\in
B_n}q^{\exc_B(\sg)}.
\end{equation}

We will establish the following relation
\begin{equation} \label{pibn}
\sum_{\pi\in B_n}q^{\exc_B(\pi)}=\sum_{k=0}^n {n \choose
k}\sum_{\sg\in D_{k}^{B}}q^{\exc_B(\sg)}.
\end{equation}
It suffices to construct a correspondence between $B_n$ and $\{(S,
T,\sigma)|S\cup T=[n]\,, S\cap T=\emptyset \, \rm{and}\, \sigma\in
B_{|S|}\}$ such that $\exc_B(\pi)=\exc_B(\sigma)$ for any $\pi\in
B_n$.

Given a signed permutation, we can decompose it into two parts
separating the fixed points from the non-fixed points. Precisely,
each $\pi \in B_n$ can be represented by $(S, T, \sg)$, where $S$ is
the set of non-fixed points of $\pi$, $T$ is the set of fixed points
of $\pi$, and $\sigma$ is a reduced signed derangement of $\pi$.
Keep in mind that $i$ is a fixed point of $\pi$ if $\pi_i=i$.  Let
$|S|=k$, then $\sigma$ is obtained  from $\pi$ as a signed
derangement on $[k]$ by deleting the fixed elements and reducing the
resulting signed permutation to $[k]$. Formally speaking, let $\tau$
be the signed permutation obtained from $\pi$ by deleting the fixed
points, then $\sigma$ is derived from $\tau$ by replacing the
minimum element by $1$ or $\bar{1}$ depending on its sign, and
replacing the second minimum element by $2$ or $\bar{2}$ depending
on its sign, and so on. Note that the elements in $\tau$ are ordered
regardless of their signs. For example, let
$\pi=4\,6\,3\,\bar{7}\,5\,1\,\bar{2}$. Then $S=\{1,2,4,6,7\}$,
$T=\{3,5\}$, $\tau=4\,6\,\bar{7}\,1\,\bar{2}$ and
$\sg=3\,4\,\bar{5}\,1\,\bar{2}$.

On the other hand, let $(S,T,\sg)$ be a representation of a signed
permutation on $[n]$ such that $S\cup T=[n]$, $S\cap T=\emptyset$,
$|S|=k$ and $\sg$ is a signed derangement on $[k]$. Then we can
recover a unique signed permutation $\pi$ on $[n]$. Let
$S=\{s_1,s_2,\ldots, s_k\}$ with $s_1<s_2<\cdots <s_k$ and $T=\{t_1,
t_2,\ldots,t_{n-k}\}$ with $t_1<t_2<\cdots<t_{n-k}$.  First $\tau$
can be obtained from $\sigma$ by replacing $i$ or $\bar{i}$ in
$\sigma$ by $s_i$ or $\bar{s_i}$, then $\pi$ can be constructed from
$\tau$ by inserting the fixed points.

We proceed to show that $\exc_B(\pi)=\exc_B(\sg)$.  Assume that $i$
is an excedance of $\pi$, that is, $\pi_i=-i$ or
$\pi_{|\pi_i|}>\pi_i$. Clearly, the fixed points do not contribute
to the number of excedances. Then there exists $j$ such that
$\pi_i=\tau_j$. If $\pi_i=-i=\tau_j$, then  there are $i-j$ fixed
points in $\pi$ before $i$. This implies that $\tau_j=-i$ is the
$j$-th minimum element in $\tau$ regardless of the signs. By the
transformation from $\tau$ to $\sigma$ as given before, we see that
$\sigma_j=-j$, that is, $j$ is an excedance of $\sigma$.

We now come to the case $\pi_{|\pi_i|}>\pi_i$. Clearly, $i$ is not a
fixed point; otherwise, we would deduce that
$\pi_{|\pi_i|}=\pi_i>\pi_i$. We have two cases depending on the sign
of $\pi_i$. First, we assume that $\pi_i$ is positive. Then $\pi_i$
is not a fixed point of $\pi$. For notational convenience, let
$\pi_i=j$, so we have $\pi_{j}>j$. Since neither $i$ nor $j$ is a
fixed point of $\pi$, both should appear in $\tau$. Thus there exist
$j_1$ and $j_2$ such that $\pi_i=\tau_{j_1}$ and $\pi_j=\tau_{j_2}$.
With the above notation, we see that $\tau_{j_2}>\tau_{j_1}$. This
implies that $\sigma_{j_2}>\sigma_{j_1}$, since the transformation
from $\tau$ to $\sigma$ is order preserving. In view of $\pi_i=j$
and $\pi_j=\tau_{j_2}$, we find that $\pi_i=j$ is the $j_2$-th
minimum element in $\tau$ regardless of the signs, namely,
$\sigma_{j_1}=j_2$. Hence we deduce that
$\sigma_{\sigma_{j_1}}=\sigma_{j_2}>\sigma_{j_1}$, namely, $j_1$ is
an excedance of $\sigma$ with $\sigma_{j_1}$ being positive.
 Conversely, given any excedance $j_1$ of
$\sigma$ with $\sigma_{j_1}$ being positive, we can reverse the
above procedure to generate an excedance $i$ of $\pi$ with $\pi_i$
being positive.

It remains to consider the case when $\pi_i=\bar{j}$  is negative.
It is clear that $i$ is not a fixed point of $\pi$. With this
notation, we have $\pi_{j}>\bar{j}$. Since $i$  is a not a fixed
point of $\pi$ and $\pi_i$ is negative, both $i$ and $\bar{j}$ will
appear in $\tau$. Hence there exist $j_1$ and $j_2$ such that
$\pi_i=\tau_{j_1}$ and $\pi_j=\tau_{j_2}$. Moreover, we see that
$\tau_{j_2}>\tau_{j_1}$ and $\sigma_{j_2}>\sigma_{j_1}$. Using the
same procedure as given before, we find that
$\sigma_{j_1}=\bar{j_2}$. It follows that
$\sigma_{|\sigma_{j_1}|}=\sigma_{j_2}>\sigma_{j_1}$, i.e., $j_1$ is
 an excedance of $\sigma$ with $\sigma_{j_1}$ being negative.
 Conversely, given an excedance ${j_1}$ of $\sigma$ with
 $\sigma_{j_1}$ being negative, we can reverse the above procedure
 to generate an excedance $i$ of $\pi$ with $\pi_i$ being negative.

Combining the above cases, we arrive at the conclusion that
$\exc_B(\pi)=\exc_B(\sg)$, which implies (\ref{pibn}). Hence we get
the following relation
\begin{equation}\label{bnq}
B_n(q)=\sum_{k=0}^n{n\choose k}d_n^B(q).
\end{equation}
Using the binomial inversion, we arrive at \eqref{d-B}. This
completes the proof. \qed

The generating function of $d_n^B(q)$ is then obtained from the
generating function of $B_n(q)$.
\begin{thm}\label{gfd}
We have
\begin{equation}\label{s2}
\sum_{n\geq
0}d_n^B(q)\frac{t^n}{n!}=\frac{(1-q)e^{tq}}{e^{2tq}-qe^{2t}}=\frac{e^{tq}}{1-\sum_{n\geq
2}2^n(q+q^2+\cdots+q^{n-1})t^n/n!}.
\end{equation}
\end{thm}

\pf Using (\ref{eg}) and (\ref{bnq}),  we obtain
\begin{equation}
e^t\sum_{n\geq 0}d_n^B(q)\frac{t^n}{n!}=\sum_{n\geq
0}B_n(q)\frac{t^n}{n!} =\frac{(1-q)e^{t(1-q)}}{1-qe^{2t(1-q)}}.
\end{equation}
The last equality of (\ref{s2}) is straightforward. This completes
the proof. \qed

Next, we give a combinatorial interpretation of the identity
\eqref{s2} based on a  generalization of the decomposition of
derangements given by Kim and Zeng \cite{Kim} for their
combinatorial proof of \eqref{gf}.

\noindent {\it A Combinatorial Proof of Theorem \ref{gfd}.} First,
we give an outline the proof of Kim and Zeng for ordinary
derangements. We adopt the convention that a cycle $\sg=s_1s_2\cdots
s_k$ of length $k$ is written in such a way that $s_1$ is the
minimum element and $\sg_{s_i}=s_{i+1}$ with $s_{k+1}=s_1$. A cycle
$\sg$ (of length at least two) is called unimodal (resp. prime) if
there exists $i$ $( 2 \leq i \leq k)$ such that $s_1 < \cdots <
s_{i-1} < s_i > s_{i+1} > \cdots> s_k$ (resp. in addition, $s_{i-1}
< s_k$). Let $(l_1,\ldots,l_m)$ be a composition of $n$, a sequence
of prime cycles $\tau= (\tau_1,\tau_2,\ldots,\tau_m)$ is called a
$P$-decomposition of type $(l_1,\ldots,l_m)$ if $\tau_i$ is of
length $l_i$ and the underlying sets of $\tau_1, \tau_2,
\ldots,\tau_m$ form a partition of $[n]$. Define the excedance of
$\tau$ as the sum of the excedances of its prime cycles, that is,
 \[ \exc (\tau)= \exc(\tau_1) + \cdots
+ \exc (\tau_m),\] and the weight of $\tau$ is defined by $q^{\exc(
\tau)}$. Note that one needs to express a cycle in the two row
permutation form for the purpose of computing the excedances. The
details are omitted here. In \cite{Kim}, Kim and Zeng obtained a
bijection which maps the number of excedances of a derangement to
the number of excedances of a $P$-decomposition with type
$(l_1,\ldots,l_m), l_i\geq 2$. Then the generating function of
$d_n(q)$ follows from the generating function of $P$-decomposition
with type $(l_1,\ldots,l_m)$, as given by
\[
{l_1+\cdots+l_m  \choose l_1,\ldots,l_m}\prod_{i=1}^m
(q+\cdots+q^{l_i-1})\frac{t^{l_1+\cdots +l_m}}{(l_1+\cdots+l_m)!}.
\]
Summing over $l_1,\ldots, l_m\geq 2$ and $m\geq 0$, we are led to
the right hand side of the generating function formula \eqref{gf}.

We now proceed to extend the above construction to type $B$
derangements. We need the cycle decomposition of a signed
permutation, which can be viewed as the cycle decomposition of an
ordinary permutation with signs attached to some elements. There is
one point that needs to be taken into account, that is, a signed
permutation is a signed derangement if and only if the cycle
decomposition does not have any one-cycle with a positive sign. More
precisely, for any derangement $\pi$ of type $B$, we can decompose
it into cycles
\[\pi=(C_1, C_2, \ldots, C_k),\]
where $C_1, C_2, \ldots, C_k$ are written in decreasing order of
their minimum elements subject to the following order
\begin{equation}\label{order}
 \bar{n}<\cdots< \bar{2}<\bar{1}<1<2<\cdots<n.
 \end{equation}
We first apply the $U$-algorithm \cite{Kim} to decompose each cycle
$\sigma=s_1s_2\cdots s_k$ of $\pi$ into a sequence of unimodal
cycles, here we impose the order (\ref{order}) in defining unimodal
cycles. It should be noted that a cycle with only one negative
element is also considered as a unimodal and prime cycle.  Then we
define $U(\pi)=(U(C_1),U(C_2),\ldots,U(C_k))$.

\noindent {\bf The $U$-algorithm}
\begin{itemize}
  \item [\rm 1.] If $\sigma$ is unimodal, then set $U(\sigma)=(\sigma)$.
  \item [\rm 2.] Otherwise, let $i$ be the largest integer such that
  $s_{i-1}>s_i<s_{i+1}$ and $j$ be the unique integer greater than
  $i$ such that $s_j>s_i>s_{j+1}$. Then set
  $U(\sigma)=(U(\sigma_1),\sigma_2)$, where $\sigma_1=s_1\cdots s_{i-1}s_{j+1}\cdots
  s_k$,  and $\sigma_2=s_is_{i+1}\cdots s_j$ is unimodal.
\end{itemize}

We claim that the number of excedances of $\pi$ is equal to total
number of excedances of unimodal cycles in $U(\pi)$.  It suffices to
verify that this statement is valid for each cycle
$\sigma=s_1s_2\cdots s_k$ of $\pi$.  Clearly, it is true if $\sigma$
is unimodal. Otherwise, it suffices to show
$\exc_B(\sigma)=\exc_B(\sigma_1)+\exc_B(\sigma_2)$. Assume $|s_t|$
is an excedance of $\sigma$, i.e.,
$\sigma_{|\sigma_{|s_t|}|}>\sigma_{|s_t|}$. Then it is necessary to
find an excedance in $\sigma_1$ or $\sigma_2$.   By the cycle
notation of $\sigma$, we have $\sigma_{|s_t|}=s_{t+1}$ and
$\sigma_{|\sigma_{|s_t|}|}=s_{t+2}$, then $|s_t|$ is an excedance of
$\sigma$ implies $s_{t+2}>s_{t+1}$. For $t=i-2$, we have
$s_{t+2}=s_i$ and $s_{t+1}=s_{i-1}$. On the other hand, we have
$s_{i}<s_{i-1}$ by the choice of $i$, so it cannot be an excedance
of $\sigma$. Using the same argument, we see that when $t=j-1$,
$|s_t|$ cannot be an excedance of $\sigma$. Therefore, if $1\leq
t<i-2$ or $j+1\leq t\leq k$, then $|s_t|$ is an excedance of
$\sigma_1$. Similarly, if $i\leq t<j-1$, then  $|s_t|$ is an
excedance of $\sigma_2$. If $t=i-1$, i.e., $|s_{i-1}|$ is an
excedance of $\sigma$, then $s_{i+1}>s_{i}$, which implies that
$|s_j|$ is an excedance of $\sigma_2$. If $t=j$, i.e., $|s_j|$ is an
excedance of $\sigma$, then $s_{j+2}>s_{j+1}$, which implies that
$|s_{i-1}|$ is an excedance of $\sigma_1$. Conversely, given an
excedance of $\sigma_1$ or $\sigma_2$, we can determine an excedance
of $\sigma$ by reversing the above procedure.

For example, let $\pi=3\,\bar{5}\,4\,2\,9\,\bar{6}\,8\,7\,\bar{1}$.
Then we have $\exc_B(\pi)=5$, and $C_1=7\,8$, $C_2=\bar{6}$ and
$C_3=\bar{5}\,9\,\bar{1}\,3\,4\,2$. Moreover, we find
\[
U(C_1)=(7\,8),\quad U(C_2)=(\bar{6}),\quad
U(C_3)=(\bar{5}\,9,\,\bar{1}\,3\,4\,2), \] and
\[U(\pi)=(7\,8,\,\bar{6},\,\bar{5}\,9,\,\bar{1}\,3\,4\,2).\]
Note that $\exc_B(U(\pi))=5$, in accordance with $\exc_B(\pi)=5$.

Next, we recall the $V$-algorithm given in \cite{Kim}, which
transforms a sequence of unimodal cycles into a sequence of prime
cycles. For signed derangements, we will use this algorithm by
imposing the order relation (\ref{order}).

\noindent {\bf The $V$-algorithm}
\begin{itemize}
   \item [\rm 1.] If $\sigma$ is prime, then set $V(\sigma)=(\sigma)$.
   \item [\rm 2.] Otherwise, let $j$ be the smallest integer such
   that $s_j>s_{i}>s_{j+1}>s_{i-1}$ for some integer $i$ greater
   than $1$. Then set $V(\sigma)=(V(\sigma_1),\sigma_2)$,
    where $\sigma_1=s_1\cdots s_{i-1}s_{j+1}\cdots
   s_k$,  and $\sigma_2=s_is_{i+1}\cdots s_j$ is prime.
\end{itemize}

We claim that the total number of excedances of  unimodal cycles in
$U(\pi)$ is equal to total number of  excedances of  prime cycles in
$V\circ U(\pi)$. It suffices to prove that the claim is valid for
any unimodal cycle $\sigma=s_1s_2\cdots s_k$. Clearly, it is true if
$\sigma$ is prime. Otherwise, without loss of generality we may show
that $\exc_B(\sigma)=\exc_B(\sigma_1)+\exc_B(\sigma_2)$. Assume
$|s_t|$ is an excedance of $\sigma$, we have $s_{t+2}>s_{t}$. As
will be shown, we can find an excedance of $\sigma_1$ or $\sigma_2$.
For $t=j-1$, $|s_t|$ cannot be an excedance of $\sigma$, since
$s_{j+1}<s_j$ by the choice of $j$. So, if $1\leq t<i-2$ or $j+1\leq
t\leq k$, then $|s_t|$ is an excedance of $\sigma_1$. If $i\leq
t<j-1$, we deduce that $|s_t|$ is an excedance of $\sigma_2$. Since
$s_i>s_{i-1}$ by the choice of $i$, when $t=i-2$, $|s_t|$ is an
excedance of $\sigma$. On the other hand, we have $s_{j+1}>s_{i-1}$
by the choices of $i$ and $j$, which implies that $|s_t|$ is an
excedance of $\sigma_1$. For $t=i-1$ or $t=j$, we can use the same
argument as in the $U$-algorithm.

Applying $V$ to each cycle of $U(\pi)$ in the above example, we
obtain that
\[
V\circ U(\pi)=(7\,8,\,\bar{6},\,\bar{5}\,9,\,\bar{1}\,2,\,3\,4),
\]
and $\exc_B(V\circ U(\pi))=5$.

 Using  the composition $V\circ U$, we can transform a derangement in $B_n$
to a $P$-decomposition of $[n]$. Moreover, it has been shown that
this map is a bijection in \cite{Kim}. In the type $B$ case, we
define the weight of each prime cycle $\tau$ by $q^{\exc_B(\tau)}$,
where $\exc_B(\tau)$ is the number of the excedances of the type $B$
derangement with cycle decomposition $\tau$. Note that in the cycle
decomposition of a type $B$ derangement, we allow cycles of length
one with negative elements. Thus the corresponding
$P$-decompositions have type $(1^k,l_1,\ldots,l_m), k\geq 0, l_i\geq
2$. For a cycle containing only one negative element, the weight is
$q$. For a cycle of length $l\geq 2$, we have $2^l$ choices for the
$l$ elements in the prime cycle, so the weight of such a prime cycle
on a  $l$-set is $2^l(q + q^2 + \cdots+ q^{l-1})$. Hence the
generating function of $d_n^B(q)$ follows from the generating
function of $P$-decompositions of type $(1^k,l_1,\ldots,l_m), k\geq
0, l_i\geq 2$, as given by
\[
q^{k}t^k{l_1+\cdots+l_m  \choose l_1,\ldots,l_m}\prod_{i=1}^m
2^{l_i}(q+\cdots+q^{l_i-1})\frac{t^{l_1+\cdots
+l_m}}{(l_1+\cdots+l_m)!}.
\]
Summing over $l_1,\ldots, l_m\geq 2$ and $k\geq 0$, $m\geq 0$, we
 obtain the right hand side of \eqref{s2}. \qed

\section{A Recurrence Relation}

In this section, we will use the recurrence relation for  Eulerian
polynomials of type $B$ to derive a recurrence relation for the
$q$-derangement polynomials $d_n^B(q)$. Applying a theorem of Zhang
\cite{Zhang2}, we deduce that $d_n^B(q)$ form a Sturm sequence, that
is, $d_n^B(q)$ has only real roots and separated by the roots of
$d_{n-1}^B(q)$. Moreover, from the initial values, one sees that
$d_n^B(q)$ has only non-positive real roots. Consequently,
$d_n^B(q)$ is log-concave. Although the polynomials $d_n^B(q)$ are
not symmetric, we show that they are almost symmetric in the sense
that the coefficients have the spiral property.

The following recurrence formula (\ref{eulerB}) for $B_n(q)$ is a
special case of Theorem 3.4 in Brenti \cite{Bre94}, see also Chow
and Gessel \cite{Chow2}, which will play a key role in establishing
a recurrence relation for $d_n^B(q)$.

\begin{thm} We have
\begin{equation}\label{eulerB}
B_n(q)=((2n-1)q+1)B_{n-1}(q)+2q(1-q){B_{n-1}'}(q), ~~~~~   n
\geqslant 1,
\end{equation}
where $B_0(q)=1$.
\end{thm}

\begin{thm}For $n\geq 2$, we have
\begin{equation}\label{dnr}
d_n^B(q)=(2n-1)qd_{n-1}^B(q)+2q(1-q){d_{n-1}^{B^\prime}}(q)+2(n-1)qd_{n-2}^B(q).
\end{equation}
\end{thm}

\pf By (\ref{d-B}) and  \eqref{eulerB}, we obtain
\allowdisplaybreaks {\footnotesize
\begin{align*}
d_n^B(q)&=
\sum_{k=0}^n (-1)^{n-k}{n \choose k}B_k(q)\\
&=\sum_{k=0}^n (-1)^{n-k}\left({n-1 \choose k-1}+{n-1 \choose
k}\right)B_k(q)\\
&=-d_{n-1}^B(q)+\sum_{k=1}^n (-1)^{n-k}{n-1 \choose k-1}(((2k-1)q+1)B_{k-1}(q)+2q(1-q)B_{k-1}'(q))\\
&=-qd_{n-1}^B(q)+2q\sum_{k=1}^n (-1)^{n-k}\left({n \choose
k}-{n-1 \choose k}\right)kB_{k-1}(q)+2q(1-q){d^{B'}_{n-1}}(q)\\
&=-d_{n-1}^B(q)+2nq\sum_{k=1}^n (-1)^{n-k}{n-1 \choose k-1}B_{k-1}(q)\\
&\qquad\, +2q(n-1)\sum_{k=1}^n(-1)^{n-k-1}{n-2 \choose k-1}B_{k-1}(q)+2q(1-q){d_{n-1}^{B'}}(q)\\
&=(2n-1)qd_{n-1}^B(q)+2(n-1)qd_{n-2}^B(q)+2q(1-q){d^{B'}_{n-1}}(q),
\end{align*}}
as desired. \qed

Equating coefficients on both sides of  \eqref{dnr}, we are led to
the following recurrence relation for the numbers $d_{n,\,k}$.

\begin{cor}
For $n\geq 2$ and $k\geq 1$, we have
\begin{equation} \label{dnkr}
d_{n,\,k}=2kd_{n-1,\,k}+(2n-2k+1)d_{n-1,\,k-1}+2(n-1)d_{n-2,\,k-1}.
\end{equation}
\end{cor}

From the above relation (\ref{dnkr}), it follows that
$d_{n,\,1}=2^n$ for $n>1$. The recurrence relation (\ref{dnr}) on
$d_n^B(q)$ enables us to show that the polynomials $d_n^B(q)$ form a
Sturm sequence. The proof turns out to be an application of the
following theorem of Zhang \cite{Zhang2}.

 \begin{thm}
 Let $f_n(q)$
be a polynomial of degree $n$ with nonnegative real coefficients
satisfying the following conditions:
\begin{itemize}
\item[\rm(1)]
For $n\geq 2$,
$f_n(q)=a_nqf_{n-1}(q)+b_nq(1+c_nq)f_{n-1}'(q)+d_nqf_{n-2}(q)$,
where $a_n>0, b_n>0, c_n\leq 0, d_n\geq 0$;
\item[\rm(2)]   For $n\geq
1$, zero is a simple root of $f_n(q)$;
 \item[\rm(3)] $f_0(q)=e,
f_1(q)=e_1q$ and $f_2(q)$ has two real roots, where $e\geq 0$ and
$e_1\geq 0$.
\end{itemize} Then the polynomial $f_n(q)$ has $n$
distinct real roots, separated by the roots of $f_{n-1}(q)$,  $n\geq
2$.
\end{thm}

It can be easily verified that the recurrence relation (\ref{dnr})
satisfies the conditions in the above theorem. Thus we reach the
following assertion.

\begin{thm}\label{roots}
The polynomials $d_n^B(q)$ form a Sturm sequence, that is,
$d_n^B(q)$ has $n$ distinct non-positive real roots, separated by
the roots of $d_{n-1}^B(q)$.
\end{thm}

As a direct consequence of the above theorem, we see that the
coefficients of $d_n^B(q)$ are log-concave. Although the
coefficients are not symmetric as in the type $A$ case, we will show
that they are almost symmetric in the sense that  they satisfy the
spiral property. The spiral property was first observed by Zhang
\cite{Zhang96} in his proof of a conjecture of Chen and Rota
\cite{Chen}.

\begin{thm}\label{spiralpro}
The polynomials $d_n^B(q)$ have the spiral property. Precisely, for
$n\geq 2$, if $n$ is even, then
\begin{equation*}
d_{n,\,n}<d_{n,\,1}<d_{n,\,n-1}<d_{n,\,2}<d_{n,\,n-2}<\cdots<d_{n,\,
\frac{n}{2}+2} <d_{n,\,\frac{n}{2}-1}<d_{n,\,\frac{n}{2}+1}<d_{n,\,
\frac{n}{2}},
\end{equation*}
and if $n$ is odd, then
\begin{equation*}
d_{n,\,n}<d_{n,\,1}<d_{n,\,n-1}<d_{n,\,2}<d_{n,\,n-2}<\cdots<d_{n,\,
\frac{n+3}{2}}<d_{n,\,\frac{n-1}{2}}<d_{n,\,\frac{n+1}{2}}.
\end{equation*}
\end{thm}

\pf Let
\[
f(n)=\begin{cases}%
\frac{n}{2}-1,&\mbox{ if $n$ is even},\\[6pt]
\frac{n-1}{2},&\mbox{ if $n$ is odd}.
\end{cases}
\]
In this notation, the spiral property can be described by the
following inequalities
\begin{equation}\label{spiral}
 d_{n,\,n+1-k}<d_{n,\,k}<d_{n,\,n-k}
\end{equation}
for any $1\leq k\leq f(n)$, and in addition, the inequality
\begin{equation} \label{n2}
d_{n,\,\frac{n}{2}+1}<d_{n,\,\frac{n}{2}}\end{equation}
 when $n$ is even.

We proceed to prove the relations (\ref{spiral}) and (\ref{n2}) by
induction on $n$. It is easily seen that  (\ref{spiral}) and
(\ref{n2}) hold for $n=2$ and $n=3$. We now assume that they hold
for all integers up to $n$. We now aim to show that
\begin{equation}\label{n+1}
d_{n+1,\,n+2-k}<d_{n+1,\,k}<d_{n+1,\,n+1-k}
\end{equation}
 for any $1\leq k\leq f(n+1)$, and it is also necessary to show that
when $n+1$ is even,
\begin{equation}\label{n21}
d_{n+1,\,\frac{n+3}{2}}<d_{n+1,\,\frac{n+1}{2}}.
\end{equation}

For $k=1$, we have $d_{n+1,\,n+1}-d_{n+1,\,1}=1-2^{n+1}<0$. For
$2\leq k\leq f(n+1)$, by the recurrence relation \eqref{dnkr} for
$d_{n,\,k}$, we have
\begin{align}
d_{n+1,\,n+2-k}&= 2(n+2-k)d_{n,\,n+2-k}+(2k-1)d_{n,\,
n+1-k}+2nd_{n-1,\,n+1-k},  \label{d1}
                                        \\[5pt]
d_{n+1,\,k}  & =2kd_{n,\,k}+(2n-2k+3)d_{n,\,k-1}+2nd_{n-1,\,k-1},\label{d2}\\[5pt]
d_{n+1,\,n+1-k}&=2(n+1-k)d_{n,\,n+1-k}+(2k+1)d_{n,\,n-k}+2nd_{n-1,\,
n-k}. \label{d3}
\end{align}
It follows from (\ref{d1}) and (\ref{d2}) that
\begin{align*}
d_{n+1,\,n+2-k}-d_{n+1,\,k}&=(2n-2k+3)(d_{n,\,n+2-k}-d_{n,\,k-1})+2k(d_{n,\,n+1-k}-d_{n,\,k})\\[5pt]
&\quad+2n(d_{n-1,\,n+1-k}-d_{n-1,\,k-1})+(d_{n,\,n+2-k}-d_{n,\,
n+1-k}).
\end{align*}
By the inductive hypothesis, we see that the difference in every
parenthesis in the above expression is negative. This implies that
for $2\leq k\leq f(n+1)$
\begin{equation}\label{nk1}
 d_{n+1,\,n+2-k}-d_{n+1,\,k}<0.
 \end{equation}
Similarly, for $2\leq k \leq f(n+1)$, in view of (\ref{d2}) and
(\ref{d3}) we find
\begin{align*}
d_{n+1,\,k}-d_{n+1,\,n+1-k}&=(2k+1)(d_{n,\,k}-d_{n,\,n-k})+2n(d_{n-1,\,k-1}-d_{n-1,\,n-k})\\[5pt]
&\quad+(2n+3-2k)(d_{n,\,k-1}-d_{n,\,n+1-k})+(d_{n,\,n+1-k}-d_{n,\,
k}).
\end{align*}
Again, by the inductive hypothesis, we deduce that for $2\leq k\leq
f(n+1)$,
\begin{equation}\label{nk2}
d_{n+1,\,k}-d_{n+1,\,n+1-k}<0.
\end{equation}
Combining (\ref{nk1}) and (\ref{nk2}) gives (\ref{n+1}) for $1\leq
k\leq f(n+1)$.

It remains to verify (\ref{n21}) when $n+1$ is even. By the
recurrence relation \eqref{dnkr} of $d_{n,k}$, we have
\begin{align*}
d_{n+1,\,\frac{n+3}{2}}&=(n+3)d_{n,\,\frac{n+3}{2}}+nd_{n,\,\frac{n+1}{2}}+2nd_{n-1,\,\frac{n+1}{2}},\\[5pt]
d_{n+1,\,\frac{n+1}{2}}&=(n+1)d_{n,\,\frac{n+1}{2}}+(n+2)d_{n,\,
\frac{n-1}{2}}+2nd_{n-1,\,\frac{n-1}{2}}.
\end{align*}
This yields
\begin{align*}
d_{n+1,\,\frac{n+3}{2}}-d_{n+1,\,\frac{n+1}{2}}&=(n+2)(d_{n,\,
\frac{n+3}{2}}-d_{n,\,\frac{n-1}{2}})+(d_{n,\,\frac{n+3}{2}}
-d_{n,\,\frac{n+1}{2}})\\[5pt]
&\qquad+2n(d_{n-1,\,\frac{n+1}{2}}-d_{n-1,\,\frac{n-1}{2}}).
\end{align*}
Again, by the inductive hypothesis, we immediately obtain
\eqref{n21}. This completes the proof.
 \qed

\section{The Limiting Distribution}

In this section, we show that the limiting distribution of the
coefficients of $d_n^B(q)$ is normal. The type $A$ case has been
studied by Clark \cite{Clar}. It has been shown that the limiting
distribution of the coefficients of $d_n(q)$ is normal. Let $\xi_n$
be the number of type $B$ excedances in a random type $B$
derangement on $[n]$. We first compute the expectation and the
variance of $\xi_n$. Then we use Lyapunov's theorem to deduce that
$\xi_n$ is asymptotically normal.

\begin{thm}
We have
\begin{align}
\E\xi_n&=\frac{n}{2}+\frac{1}{4}+o(1),\\
 \Var\xi_n&=\frac{n}{12}-\frac{1}{16}+o(1).
\end{align}
\end{thm}

\pf By the recurrence relation \eqref{eulerB} for $B_n(x)$, we have
for $n\geq 1$,
\begin{equation}\label{bn1}
B_n'(x)=(2n-1)B_{n-1}(x)+(2nx-5x+3)B_{n-1}'(x)+2x(1-x)B_{n-1}''(x).
\end{equation}
Since $B_n(1)=2^nn!$ for $n\geq 0$, setting $x=1$ in \eqref{bn1}
gives the following recurrence relation for $B_n'(1)$:
\begin{equation*}
B_n'(1)=(2n-1)(n-1)!2^{n-1}+(2n-2)B_{n-1}'(1).
\end{equation*}
It can be verified that for $n\geq 1$,
\begin{equation}\label{b'1}
\ B_n'(1)=\frac{n2^nn!}{2}.
\end{equation}
Moreover, by (\ref{bn1}) we get
\begin{equation}\label{bn2}
B_n''(x)=(4n-6)B_{n-1}'(x)+(2nx-9x+5)B_{n-1}''(x)+2x(1-x)B_{n-1}'''(x).
\end{equation}
Setting  $x=1$ in \eqref{bn2} and using  \eqref{b'1}, we obtain
\begin{equation}
B_n''(1)=(2n-3)(n-1)2^{n-1}(n-1)!+(2n-4)B_{n-1}''(1).
\end{equation}
Then for $n \geq 2$,
\begin{equation}\label{b''1}
\ B_n''(1)=\frac{(3n^2-5n+1)2^nn!}{12}.
\end{equation}
Since $B_n(1)=2^nn!$,  in view of the formula \eqref{d-B}, we see
that
\begin{equation}\label{dnb1}
d_n^B(1)=\sum_{k=0}^n(-1)^{n-k}{n\choose k}B_k(1).
\end{equation}
Let
 \[ s_n=\sum_{k=0}^n(-1)^k\frac{1}{2^k\cdot k!}.\]
 Then $d_n^B(1)$ can be written as $2^nn!s_n$.

Applying the formula \eqref{d-B} again and using the evaluation
\eqref{b'1} for $B_n'(1)$, we find that for $n\geq 1$
\begin{align}\nonumber
{d_n^B}^\prime(1)&=\sum_{m=1}^n(-1)^{n-m}{n\choose m}\cdot \nonumber
2^{m-1}m\cdot m!\\ \nonumber
&=2^nn!\sum_{m=0}^{n-1}(-1)^m\frac{n-m}{m!2^{m+1}}\\ \nonumber
&=2^nn!\left(\frac{n}{2}\sum_{m=0}^{n-1}(-1)^m\frac{1}{m!2^m}+\frac{1}{4}\sum_{m=0}^{n-1}(-1)^{m-1}\frac{1}{2^{m-1}(m-1)!}\right)\\
&=\frac{2^nn!}{2}\left(ns_{n-1}+\frac{1}{2}s_{n-2}\right).\label{dnb11}
\end{align}

Differentiating \eqref{d-B} twice  and invoking \eqref{b''1}, we
deduce that for $n\geq 2$
\begin{align}\nonumber
{d_n^B}^{\prime\prime}(1)= &\sum_{m=2}^n(-1)^{n-m}{n\choose m}2^m
m!\frac{3m^2-5m+1}{12}\\ \nonumber
=&\frac{2^nn!}{12}\sum_{m=2}^n(-1)^{n-m}\frac{3m^2-5m+1}{(n-m)!2^{n-m}}\\
\nonumber
=&\frac{2^nn!}{12}\sum_{m=0}^{n-2}(-1)^m\frac{3(n-m)^2-5(n-m)+1}{m!2^m}\\
\nonumber
=&\frac{2^nn!}{12}\left(\sum_{m=0}^{n-2}(-1)^m\frac{3n^2-5^n+1}{m!2^m}+\sum_{m=0}^{n-2}(-1)^m\frac{-6n+5}{(m-1)!2^m}\right.
\\ \nonumber
&\qquad\quad\left.+3\sum_{m=0}^{n-2}(-1)^m\frac{m^2}{m!2^m}\right)\\
\nonumber
=&\frac{2^nn!}{12}\left((3n^2-5n+1)s_{n-2}+\frac{1}{2}(6n-5)s_{n-3}+
\frac{3}{4}s_{n-4}-\frac{3}{2}s_{n-3}\right)\\
=&\frac{2^nn!}{12}\left((3n^2-5n+1)s_{n-2}+(3n-4)s_{n-3}+\frac{3}{4}s_{n-4}\right).\label{dnb12}
\end{align}

It is easy to see that $s_{n-r}/s_n=1+o(1)$ for  $r=1,2,3,4$. From
\eqref{dnb1}, \eqref{dnb11} and \eqref{dnb12}, we conclude that
\begin{align}\label{exp}
&\E\xi_n=\frac{{d_n^B}'(1)}{d_n^B(1)}=\frac{n}{2}+\frac{1}{4}+o(1),\\
\label{var}
&\Var\xi_n=\frac{{d_n^B}''(1)}{d_n^B(1)}+\E\xi_n-(\E\xi_n)^2=\frac{n}{12}-\frac{1}{16}+o(1),
\end{align}
as desired. \qed

Given the formulas for the expectation and variance of $\xi_n$, we
will use
 Lyapunov's theorem \cite[Section 1.2]{sach}
 to derive that the limiting distribution of $\xi_n$ is normal.
Recall that a triangular array of independent random variables
$\xi_{nk},\, k=1,2,\ldots,n$, $n=1,2,\ldots$, is called a Poisson
sequence if
\begin{equation*}
P\{\xi_{nk}=1\}=p_k,\qquad P\{\xi_{nk}=0\}=q_k,
\end{equation*}
where $p_k=p_k(n)$, $q_k=q_k(n)$ and $p_k+q_k=1$. Then Lyapunov's
theorem can be used to derive asymptotically normal distributions
\cite[Section 1.2]{sach}.

\begin{thm}[Lyapunov]\label{lya} Let
\begin{equation*}
V_n^2=\sum_{k=1}^np_kq_k,\qquad
\eta_n=V_n^{-1}\sum_{k=1}^n(\xi_{nk}-p_k).
\end{equation*}
If $V_n\rightarrow \infty$ as $n\rightarrow\infty$, then the
sequence $\{\eta_n\}$ is asymptotically standard normal.
\end{thm}

The above theorem enables us to derive the asymptotic distribution
of the random variable $\eta_n$.

\begin{thm}\label{normaldis}
The distribution of the random variable
\begin{equation*}
\eta_n=\frac{\xi_n-\E\xi_n}{\sqrt{\Var\xi_n}}
\end{equation*}
converges to the standard normal distribution as
$n\rightarrow\infty$.
\end{thm}

\pf Since the polynomials $d_n^B(q)$ have distinct, real and
non-positive roots, we may express $d_n^B(q)$ as
\[ d_n^B(q)=q(q+\alpha_1)
(q+\alpha_2)\cdots (q+\alpha_{n-1}).\] Hence the random variable
$\xi_n$, namely, the number of type $B$ excedances in a random type
$B$ derangement on $[n]$, can be represented as the sum of
independent random variables
\begin{equation*}
\xi_n=\xi_{n1}+\xi_{n2}+\cdots+\xi_{n,n-1}+\xi_{n,n},
\end{equation*}
where $\xi_{n1}, \xi_{n2}, \ldots, \xi_{n,n-1},\xi_{n,n}$ form a
Poisson sequence with $p_k=\frac{1}{1+\alpha_k}$, for $k=1,2,\ldots,
n-1$ and $p_n=1$. On the other hand, the formula \eqref{var} implies
that $\Var\xi_n\rightarrow \infty$ as $n\rightarrow \infty$.
Therefore, from  Theorem \ref{lya} we deduce that $\eta_n$ is
asymptotically standard normal. \qed

\vspace{.2cm} \noindent{\bf Acknowledgments.} This work was
supported by  the 973 Project, the PCSIRT Project of the Ministry of
Education, the Ministry of Science and Technology, and the National
Science Foundation of China.


\end{document}